\documentclass[12pt,intlim,righttag]{article}
\usepackage{amsmath,amsthm}
\usepackage{amssymb}
\usepackage{amsfonts}

\newcommand{\esssup}{\operatornamewithlimits{ess\,sup}}

\newcommand{\la}{\langle}
\newcommand{\ra}{\rangle}

\newcommand{\gm}{\gamma}

\newcommand{\td}{\tilde}

\newcommand{\be}{\begin}
\newcommand{\ee}{\end}
\newcommand{\lbl}{\label}

\newcommand\beq{\begin{equation}}
\newcommand\eeq{\end{equation}}

\newcommand{\bea}{\begin{eqnarray}}
\newcommand{\eea}{\end{eqnarray}}
\newcommand{\beaa}{\begin{eqnarray*}}
\newcommand{\eeaa}{\end{eqnarray*}}

\theoremstyle{Theorem}

\theoremstyle{corollary}

\theoremstyle{remark}

\theoremstyle{definition}

\def\a{\alpha}
\def\b{\beta}
\def\g{\gamma}

\def\z{\zeta}

\def\l{\lambda}

\def\t{\tau}

\def\o{\omega}

\def\t{\tau}

\def\L{\Lambda}

\def\cE{{\cal E}}
\def\cF{{\cal F}}

\def\cM{{\cal M}}

\begin{document}
\title{The Adomian series representation of some quadratic BSDEs
 }

\author{R. Tevzadze}

\date{~}
\maketitle

\be{abstract}
{\bf Abstract.}{ The representation of the solution of some Backward Stochastic Differential Equation as an infinite series is obtained.
Some exactly solvable examples are considered.}
\ee{abstract}


\


\

\section{\bf Introduction}

In a number of papers[1,2] Adomian develops a numerical technique using
special kinds of polynomials for solving non-linear functional equations.
However, Adomian and his collaborators did not develop widely the problem of convergence.

In this article we will study by Adomian technique some kind of quadratic backward martingale equation and prove the convergence of the series. For example we tackle   an equation of the form
\beq\lbl{01}
{\cal E}_T(m){\cal E}^\alpha_T(m^\bot)=c\exp\{\eta\}
\eeq
w.r.t. stochastic integrals $m=\int f_sdW_s,\;m^\bot=\int g_sdW_s^\bot$ and real number $c$, where $(W,W^\bot)$ is 2-dimension Brownian Motion and $\eta$ is a random variable.

Equations of such type are arising in mathematical ﬁnance and they are used to characterize
optimal martingale measures (see, Biaginiat at al (2000), Mania and Tevzadze
(2000), (2003),(2006)). Note that equation (1) can be applied also to the ﬁnancial market
models with inﬁnitely many assets (see M. De Donno  at al (2003)). In
Biagini at al (2000) an exponential equation of the form
$$
\frac{\cE_T(m)}{\cE_T(m^\bot)}=ce^{\int_0^T\l^2_sds}
$$
was considered (which corresponds to the case $\alpha=-1$ ).

Our goal is to show the solvability of the equation (\ref{01}) using the Adomian method  proving the convergence of series. On the one hand, a simpler proof of solvability is obtained.
On the other hand, it allows to obtain the approximation of  the solution. It is possible to find a solution in the form of series, if we define a sequence of martingales w.r.t. the measure $\cE_T(\sum_i^nm_i+\sum_i^nm_i^\bot)\cdot P$
from  equations $c'\cE_T(m_{n+1}'+m_{n+1}'^\bot)=\cE_T^2(m_n'^\bot)$, where $m_{n+1}'= m_{n+1}-\la m_{n+1},\sum_i^n,m_{i}\ra,$ $ m_{n+1}'^\bot=m_{n+1}^\bot -\la m_{n+1}^\bot,\sum_i^nm_{i}^\bot\ra$,
and then we write down the solution
$$m=\sum_k^\infty m_k,\;m^\bot=\sum_k^\infty m_k^\bot$$ provided  the series are convergent. The proof of the convergence  is greatly simplified if we present equation
 as a BSDE in the space of BMO-martingales and use the properties of the BMO-norm. The result is resumed in Theorem 1.

 Finally we provide some examples, exactly solvable by Adomian series and also example non-solvable at all.

\section{The main result}

 Let $(\Omega,{\cal F},P)$ be a probability space with filtration
${\bf F}=({\cal F}_t,t\in[0,T])$. We assume that all local
martingales with respect to ${\bf F}$ are continuous. Here $T$ is
a fixed time horizon and ${\cal F}={\cal F}_T$.

 Let ${\cal M}$ be a stable subspace of the  space of square
 integrable martingales $H^2$. Then its ordinary orthogonal ${\cal M}^\bot$
 is a stable subspace and any element of ${\cal M}$ is strongly
 orthogonal to any element of ${\cal M}^\bot$ (see, e.g. \cite{DM}, \cite{J}).

  We consider the following exponential equation
 \begin{equation}\label{i1}
 {\cal E}_T(m){\cal E}^\alpha_T(m^\bot)=c\exp\{\eta\},
 \end{equation}
where $\eta$ is a given $F_T$-measurable random variable and $\alpha$ is a given real
number. A solution of equation (\ref{i1}) is a triple $(c,m,m^\bot)$, where  $c$ is
strictly positive constant, $m\in\cal M$ and $m^\bot\in\cal M^\bot$.
Here ${\cal E}(X)$ is the Doleans-Dade exponential of $X$.

It is evident that if $\alpha=1$ then equation (\ref{i1}) admits an "explicit" solution.
E.g., if $\alpha=1$ and $\eta$ is bounded, then using the unique
decomposition of the martingale $E(\exp\{\eta\}/F_t)$
\begin{equation}\label{i0}
E(\exp\{\eta\}/F_t)=E\exp\{\eta\}+m_t(\eta)+ m_t^\bot(\eta),\;\;m(\eta)\in{\cal M},\;\;
m^\bot(\eta)\in{\cal M}^\bot,
\end{equation}
it is easy to verify that the triple $c=\frac{1}{E\exp\{\eta\}}$,
$$
m_t=\int_0^t\frac{1}{E(\exp\{\eta\}/F_s)}dm_s(\eta),\;\;
m_t^\bot=\int_0^t\frac{1}{E(\exp\{\eta\}/F_s)}dm^\bot_s(\eta)
$$
satisfies equation  (\ref{i1}).

Our aim is to prove the existence of a unique solution of equation (\ref{i1}) for
arbitrary $\alpha\neq 0$ and $\eta$
of a general structure, assuming that it satisfies the following boundedness condition:

B) $\eta$ is an $F_T$-measurable random variable of the form
\begin{equation}\label{B}
\eta= \bar\eta+ \gamma A_T,
\end{equation}
where $\bar\eta\in L^\infty$, $\gm$ is a constant
and $A=(A_t, t\in[0,T])$ is a continuous $F$-adapted
 process of finite variation such that
$$
E(var_T(A)-var_\tau(A)/F_\tau)\le C
$$
for all stopping times $\tau$ for a constant $C>0$.

One can show that equation (\ref{i1}) is equivalent to the following semimartingale backward equation
with the square generator
\begin{eqnarray} \label{ip5}
Y_t=Y_0-\frac{\gm}{2} A_t-\la L\ra_t-\frac1\a \la L^\bot\ra_t+L_t+L^\bot_t,\;\;\;\;
Y_T=\frac{1}{2}\bar\eta.
\end{eqnarray}

We use also the equivalent equation of the form
\[
L_T+L^\bot_T=c+\la L\ra_T+\frac1\a \la L^\bot\ra_T+\frac{\gm}{2} A_T.
\]
w.r.t. $(c,L,L^\bot)$.

We use notations $|M|_{{}_{\rm BMO}}=\inf\{C:E^\frac12(\la M\ra_T-\la M\ra_\t|\cF_\t)\le C\}$ for BMO-norms of martingales,  $|A|_\o=\inf\{C:E(var_t^T(A)|\cF_t)\le C\}$ for norms of finite variation processes
and $A\cdot M$ for stochastic integrals.

Let us consider the system of semimartingale backward equations
\beaa
Y_t^{(0)}=Y_0^{(0)}-\frac{\gm}{2} A_t+L_t^{(0)}+L_t^{(0)\bot},\;\;\;
Y_T^{(0)}=\frac{1}{2}\bar\eta,\\
Y_t^{(n+1)}=Y_0^{(n+1)}\\
-\sum_{k=0}^n\la L^{(k)},L^{(n-k)}\ra_t-\frac1\a \sum_{k=0}^n\la L^{(k)\bot},L^{(n-k)\bot}\ra_t+L_t^{(n+1)}+L_t^{(n+1)\bot},\\
Y_T^{(n+1)}=0.
\eeaa
The sequence $Y_0^{(n)}=c^{(n)}, L^{(n)}+L^{\bot(n)}, n=0,1,2,\cdots$ can be defined consequently by the equations
\beaa
E(\eta|\cF_t)+\frac{\gm}{2} E(A_T|\cF_t)=c^{(0)}+L_t^{(0)}+L_t^{\bot(0)},\\
\sum_{k=0}^nE(\la L^{(k)},L^{(n-k)}\ra_T|\cF_t)-\frac1\a \sum_{k=0}^nE(\la L^{(k)\bot},L^{(n-k)\bot}\ra_T|\cF_t)\\=c^{(n+1)}
+L_t^{(n+1)}+L_t^{\bot(n+1)}.
\eeaa
{\bf Remark.} If $A_t=\int_0^ta(s,W_s,B_s)ds$, then the solution of (\ref{ip5}) is of the form $Y_t=v(t,W_t,B_t)$, where  $v(t,x,y)$ is decomposed as series $\sum_nv^n(t,x,y)$ satisfying the system of PDEs
\beaa
(\partial_t+\frac12\Delta)v^0(t,x,y)+a(t,x,y)=0,\;\;\;\;v^0(T,x,y)=0,\\
(\partial_t+\frac12\Delta)v^n(t,x,y)\\
+\frac12\sum_{k=0}^{n-1}(v_x^k(t,x,y)v_x^{n-k-1}(t,x,y)+\a v_y^k(t,x,y)v_y^{n-k-1}(t,x,y))=0,\\
v^n(T,x,y)=0,\;n\ge1.
\eeaa

\be{lem} Let
$$
Y_t=Y_0+A_t+m_t,\;\;\;\;Y_T=\eta,
$$
where $ m$ is a martingale, $\eta\in L_\infty$ and $|A|_\omega<\infty$. Then $m\in BMO$ and
\begin{equation}\label{est1}
|m|_{_{\rm BMO}}\le |\eta |_\infty+|A|_\omega.
\end{equation}
In particular,  if $|A|_\omega <\infty$ then the martingale $E(A_T|F_t)$ belongs to the $BMO$ space and
$$
|E(A_T|F_.)|_{_{\rm BMO}}\le |A|_\omega.
$$
\ee{lem}
{\it Proof.} By the Ito formula
$$
Y_t^2= 2\int_0^tY_s dm_s+2\int_0^tY_sdA_s+\la m\ra_t.
$$
Taking the difference $Y^2_\tau-Y^2_T$ and conditional expectations we have that
$$
Y^2_\tau+E(\la m\ra_{T}-\la m\ra_{\tau}|F_\tau)=E(\eta^2|F_\tau)-2E(\int_\tau^TY_sdA_s|F_\tau)\le
$$
\begin{equation}\label{est2}
\le |\eta |_\infty^2+2|Y|_\infty |A|_\omega.
\end{equation}
$E(\int_\tau^TY_sdm_s|F_\tau)=0$, since $Y_t\le E(\eta+|A_T-A_t||\cF_t)$ is bounded and $m$ is a martingale.
 Since the right-hand side of (\ref{est2}) does not depend on $\tau$ from (\ref{est2}) we obtain
$$
|Y|^2_\infty+||m||_{BMO}^2\le  |\eta |_\infty^2+|Y|^2_\infty+ |A|^2_\omega.
$$
Therefore
$$
||m||_{BMO}^2\le  |\eta |_\infty^2+ |A|^2_\omega,
$$
which implies inequality (\ref{est1}).

\be{lem} For the $BMO$ norms of martingales $L^{(n)}+L^{\bot(n)}$, defined above, the following estimates
are true
\begin{equation}\label{2}
|L^{(n)}+L^{\bot(n)}|_{_{\rm BMO}}\le a_n(1+|\b|)^{n}|L^{(0)}+L^{\bot(0)}|_{{}_{\rm BMO}}^{n+1},
\end{equation}
where the coefficients $a_n$ are calculating recurrently from
$$
a_0=1,\;\;\; a_{n+1}=\sum_{k=0}^n a_k a_{n-k}.
$$
\ee{lem}
{\it Proof.} Using Lemma 1 it is easy to show that
$$
|L^{(1)}+L^{\bot(1)}|_{_{\rm BMO}}\le a_1(1+|\b|)|L^{(0)}+L^{\bot(0)}|_{{}_{\rm BMO}}^{2},
$$
$$
|L^{(2)}+L^{\bot(2)}|_{_{\rm BMO}}\le a_2(1+|\b|)^2|L^{(0)}+L^{\bot(0)}|_{{}_{\rm BMO}}^{3}.
$$
Assume that inequality (\ref{2}) is valid for any $k\le n$ and let us show that
\begin{equation}\label{3}
|L^{(n+1)}+L^{\bot(n+1)}|_{_{\rm BMO}}\le a_{n+1}(1+|\b|)^{n+1}|L^{(0)}+L^{\bot(0)}|_{{}_{\rm BMO}}^{n+2}.
\end{equation}
Applying  Lemma 1 for $Y_t^{(n+1)}$ and the Kunita-Watanabe inequality we have
$$
|L^{(n+1)}+L^{\bot(n+1)}|_{_{\rm BMO}}\le
$$
$$
\le\esssup_\t \sum_{k=0}^n E(var_\t^T(\sum_k^n\la L^{(k)},L^{(n-k)}\ra+\b \la L^{\bot(k)},L^{\bot(n-k)}\ra)|\cF_\t)
$$
$$
\le \sum_{k=0}^n \esssup_\t E^\frac12(var_\t^T\la L^{(k)}\ra|\cF_\t)E^\frac12(var_\t^T\la L^{\bot(n-k)}\ra|\cF_\t)
$$
$$
+|\b|\sum_{k=0}^n \esssup_\t E^\frac12(var_\t^T\la L^{\bot(k)}\ra|\cF_\t)E^\frac12(var_\t^T\la L^{\bot(n-k)}\ra|\cF_\t)
$$
$$
\le\sum_k^n |L^{(k)}|_{{}_{\rm BMO}}|L^{(n-k)}|_{{}_{\rm BMO}}+|\b||L^{\bot(k)}|_{{}_{\rm BMO}}|L^{\bot(n-k)}|_{{}_{\rm BMO}}
$$
\begin{equation}\label{4}
\le(1+|\b|)\sum_{k=0}^n |L^{(k)}+L^{\bot(k)}|_{{}_{\rm BMO}}|L^{(n-k)}+L^{\bot(n-k)}|_{{}_{\rm BMO}}.
\end{equation}
Therefore, from (\ref{4}), using inequalities (\ref{2})  for any $k\le n$,  we obtain
$$
|L^{(n+1)}+L^{\bot(n+1)}|_{_{\rm BMO}}\le
$$
$$
\le(1+|\b|)\sum_{k=0}^n a_k (1+|\b|)^k|L^{(0)}+L^{\bot(0)}|^{k+1}_{{}_{\rm BMO}}a_{n-k} (1+|\b|)^{n-k}||L^{(n-k)}+L^{\bot(n-k)}|^{n-k+1}_{{}_{\rm BMO}}
$$
$$
\le(1+|\b|)^{n+1}|L^{(0)}+L^{\bot(0)}|^{n+2}_{{}_{\rm BMO}}\sum_{k=0}^n a_k a_{n-k}=
$$
$$
=a_{n+1}(1+|\b|)^{n+1}|L^{(0)}+L^{\bot(0)}|^{n+2}_{{}_{\rm BMO}}
$$
and the validity of inequality (\ref{2})  follows by induction.

\be{thr} The series $\sum_{n\ge0}(L^{(n)}+L^{\bot(n)})$ is convergent in BMO-space, if $\g$ and $|\bar\eta|_\infty$ are small enough and the sum of series is a solution of the equation (\ref{ip5}).
\ee{thr}
{\it Proof.} Without loss of generality assume that $\eta=0$.
Using the lemma 2  we get
\[
|L^{(n)}+L^{\bot(n)}|_{_{\rm BMO}}\le a_n(1+|\b|)^{n}|L^{(0)}+L^{\bot(0)}|_{{}_{\rm BMO}}^{n+1}\le a_n(1+|\b|)^{n}|\g A|_\o^{n+1}.
\]
By  lemma 3 of appendix, since $$\overline\lim_{n\to\infty}\sqrt[n]a_n=\overline\lim_{n\to\infty}\sqrt[n]{\frac1{2n+1}C^{2n+2}_{n+1}}=\overline\lim_{n\to\infty}\sqrt[n]{\frac{(2n)!}{n!n!}}
=\overline\lim_{n\to\infty}\sqrt[n]{\frac{(2n)^{2n}}{n^{2n}}}=4,$$
 the series   is convergent, when $\g<\frac{1}{4|A|_\o(1+|\b|)} $.

{\bf Remark.} Since
$\max(|L|_{{}_{\rm BMO}},|L^{\bot}|_{{}_{\rm BMO}})\le|L+L^{\bot}|_{{}_{\rm BMO}}\le |L|_{{}_{\rm BMO}}+|L^{\bot}|_{{}_{\rm BMO}}$
the convergence $\sum_{n\ge0}(L^{(n)}+L^{\bot(n)})$ implies convergence of $\sum_{n\ge0}L^{(n)}$ and $\sum_{n\ge0}L^{\bot(n)}$ and vice versa.

The existence of the solution for arbitrary bounded $\eta$ is proven  \cite{rm}.  We can prove here little more general result

\be{prop} There exists solution of (\ref{i1}) for sufficiently small $\g$ and arbitrary bounded $\bar\eta$ .
\ee{prop}
{\it Proof.} Let $\bar m+\bar m^\bot$ be solution of (\ref{i1}) for $\eta=\gm A_T$ and sufficiently small $\g$. From the result of  \cite{rm}
there exists a solution of
\[
{\cal E}_T(\td m){\cal E}^\alpha_T(\td m^\bot)=c\exp\{\bar\eta\},
\]
w.r.t $$\bar P=\cE_T(\bar m+\bar m^\bot)\cdot,\;\td m+\td m^\bot\in \cM(F,\bar P)+\cM^\bot(F,\bar P).P.$$
It is easy to verify that $m+m^\bot=\bar m+\bar m^\bot+\td m+\td m^\bot$ is a solution of (\ref{i1}) for $\eta=\bar\eta+\gm A_T$.

The uniqueness of the solution was proved in \cite{rm}.

\be{prop}. Let $\eta$ be an $\cF_T$ -measurable random variable. If there exists a triple
$(c, m, m^\bot)$, where $c \in R_+, m \in BMO \cap \cM, m^\bot \in BMO \cap \cM^\bot$ satisfying equation (\ref{i1})
then such solution is unique.
\ee{prop}
We now show that without finiteness of $|A|_\o$ either the solution does not exists or the convergence of series is valid in a week sense.

{\bf Example 1.} Let $\a=-1,\;\gm=2,\;\bar\eta=0,\;A_t=\frac12\int_0^t(W_s^2+W_s^{2\bot})ds,\;\;{\bf F}=({\cal F}_t^{W,W^\bot}),$ where $W,W^\bot$ is 2-dimensional Brownian motion. Then (\ref{ip5}) becomes
\[
L_T+L_T^\bot=c+\la L\ra_T-\la L^\bot\ra_T+\frac12\int_0^T(W_s^2+W_s^{2\bot}) ds.
\]
We have
\beaa
L_T^{(0)}+L_T^{(0)\bot}=c_0+\int_0^T(T-s)W_s dW_s+\int_0^T(T-s)W_s^\bot dW_s^\bot,\\
L_T^{n+1)}+L_T^{(n+1)\bot}=c_n+\sum_{k=0}^n\la L^{(k)},L^{(n-k)}\ra_T-\sum_{k=0}^n\la L^{(k)\bot},L^{(n-k)\bot}\ra_T,\;n\ge0.\\
\eeaa
Let assume
\beaa
L_T^{(n)}=\int_0^T(T-s)^{2n+1}\a_nW_tdW_s,\\
L_T^{(n)\bot}=\int_0^T(T-s)^{2n+1}\b_nW_t^\bot dW_s^\bot.
\eeaa
Then $a_0=1,\;\;\b_0=1$ and
\beaa
L_T^{(n+1)}=c_n'
+\sum_{k=0}^n\int_0^T(T-s)^{2n+2}\a_k\a_{n-k}W_s^2ds\\
L_T^{(n+1)\bot}=c_n''-\sum_{k=0}^n\int_0^T(T-s)^{2n+2}\b_k\b_{n-k}W_s^2ds,\;n\ge0.
\eeaa
Taking stochastic derivatives $D_t,D_t^\bot$ and conditional expectations on both sides we get
\beaa
(T-s)^{2n+3}\a_nW_t=2\sum_{k=0}^n\a_k\a_{n-k}W_t\int_t^T(T-s)^{2n+2}ds\\=\frac2{2n+3}W_t(T-t)^{2n+3}\sum_{k=0}^n\a_k\a_{n-k},\\
(T-s)^{2n+3}\b_nW_t^\bot=-\frac2{2n+3}W_t^\bot(T-t)^{2n+3}\sum_{k=0}^n\b_k\b_{n-k},
\eeaa
which means that
\[
\a_{n+1}=\frac2{2n+3}\sum_{k=0}^n\a_k\a_{n-k},\; \b_{n+1}=-\frac2{2n+3}\sum_{k=0}^n\b_k\b_{n-k},n\ge0.
\]
Introducing $\a(s)=\sum_{n=0}^\infty\a_n  s^{2n+1},\;\b(s)=\sum_{n=0}^\infty\b_n  s^{2n+1}$ one obtains
\beaa
\a'(s)=\a_0+\sum_{n=0}^\infty(2n+3)\a_{n+1}s^{2n+2}\\
=1+2\sum_{n=0}^\infty\sum_{k=0}^n(\a_k\a_{n-k})s^{2n+2}=1+2a^2(s),\\
\b'(s)=\b_0+\sum_{n=0}^\infty(2n+3)\b_{n+1}  s^{2n+2}\\
=1-2\sum_{n=0}^\infty\sum_{k=0}^n\b_k\b_{n-k}  s^{2n+2}=1-2\b^2(s).
\eeaa
I.e.
\bea\lbl{ode}
\a'(s)=1+2a^2(s),\;\a(0)=0,\\
\notag \b'(s)=1-2\b^2(s),\;\b(0)=0.
\eea
Thus
\[
\a(s)=\frac1{\sqrt2}\tan(\sqrt2 s),\;\b(s)=-\frac1{\sqrt2}\tanh(\sqrt2 s).
\]
If  $T<\frac{\pi}{2\sqrt2}$ series are convergent (not in BMO-space) and $(c,L,L^\bot)$ is defined as $c=\frac12\ln {\cos(\sqrt2 T)}{\cosh(\sqrt2 T)}$ (by calculations in the appendix),
$$L_t=\frac1{\sqrt2}\int_0^t\tan(\sqrt2 s)W_sdW_s,\;L_t^\bot=-\frac1{\sqrt2}\int_0^t\tanh(\sqrt2 s)W_s^\bot W_s^\bot.$$
When $T>\frac{\pi}{2\sqrt2}$ a local martingale $L$ satisfying  $L_T-\la L\ra_T=\frac12\int_0^TW_t^2dt$ does not exist (despite the fact that $\int_0^TW_t^2dt$ is p-integrable for each $p\ge1$),
since from $\cE_T(2L)=e^{\int_0^TW_t^2dt}$ follows that
$Ee^{\int_0^TW_t^2dt}=E\cE_T(2L)\le1$, which contradicts to $Ee^{\int_0^TW_t^2dt}=\infty$ (see appendix).

In the next example exact solution of (\ref{ip5}) also exists, however it does not belong  to the extreme cases considered in \cite{rm0},\cite{rm1}.

{\bf Example 2.} Let $\a=-1,\;\gm=2,\;\bar\eta=0,\;A_t=\int_0^tW_sW_s^\bot ds,\;\;{\bf F}=({\cal F}_t^{W,W^\bot}),$ where $W,W^\bot$ is a 2-dimensional Brownian motion. Then (\ref{ip5}) becomes
\[
L_T+L_T^\bot=c+\la L\ra_T-\la L^\bot\ra_T+\int_0^TW_sW_s^\bot ds.
\]
We have
\beaa
L_T^{(0)}=EL_T^{(0)}+\int_0^T(T-s)W_s^\bot dW_s,\;L_T^{(0),\bot}=EL_T^{(0),\bot}+\int_0^T(T-s)W_s dW_s^\bot,\\
L_T^{(n+1)}+L_T^{(n+1)\bot}=c_n+\sum_{k=0}^n\la L^{(k)},L^{(n-k)}\ra_T-\sum_{k=0}^n\la L^{(k)\bot},L^{(n-k)\bot}\ra_T,\;n\ge0.\\
\eeaa
We assert that
\beaa
L_T^{(n)}=EL_T^{(n)}+\int_0^T(T-s)^{2n+1}(\a_nW_t+\b_nW_s^\bot )dW_s,\\
L_T^{(n)\bot}=EL_T^{(n)\bot}+\int_0^T(T-s)^{2n+1}(\b_nW_t-\a_nW_s^\bot )dW_s^\bot,
\eeaa
where $\a_0=0,\;\;\b_0=1$ and
\[
\a_{n+1}=\frac2{2n+3}\sum_{k=0}^n(\a_k\a_{n-k}-\b_k\b_{n-k}),\; \b_{n+1}=\frac4{2n+3}\sum_{k=0}^n\a_k\b_{n-k},n\ge0.
\]
Indeed,
\beaa
L_T^{(n+1)}+L_T^{(n+1)\bot}=c_n\\
+\sum_{k=0}^n\int_0^T(T-s)^{2n+2}(\a_kW_s+\b_kW_s^\bot )(\a_{n-k}W_s+\b_{n-k}W_s^\bot )ds\\
-\sum_{k=0}^n\int_0^T(T-s)^{2n+2}(\b_kW_s-\a_kW_s^\bot )(\b_{n-k}W_s-\a_{n-k}W_s^\bot )ds\\
=\sum_{k=0}^n\int_0^T(T-s)^{2n+2}[(\a_k\a_{n-k}-\b_k\b_{n-k})W_s^2-(\a_k\a_{n-k}-\b_k\b_{n-k})W_s^{\bot2}\\
+2(\a_k\b_{n-k}+\b_k\a_{n-k})W_sW_s^\bot]ds+c_n,\;n\ge0.
\eeaa
Using representation of integrands by stochastic derivatives we get
\beaa
(T-t)^{2n+3}(\a_{n+1}W_t+\b_{n+1}W_t^\bot )\\
=E[D_t(\sum_{k=0}^n\la L^{(k)},L^{(n-k)}\ra_T-\sum_{k=0}^n\la L^{(k)\bot},L^{(n-k)\bot}\ra_T)|\cF_t]\\
=2\sum_{k=0}^n[(\a_k\a_{n-k}-\b_k\b_{n-k})W_t+(\a_k\b_{n-k}+\b_k\a_{n-k})W_t^\bot]\int_0^T(T-s)^{2n+2}ds\\
=\frac{2(T-t)^{2n+3}}{2n+3}\sum_{k=0}^n[(\a_k\a_{n-k}-\b_k\b_{n-k})W_t+(\a_k\b_{n-k}+\b_k\a_{n-k})W_t^\bot],\\
(T-t)^{2n+3}(\b_{n+1}W_t-\a_{n+1}W_t^\bot)\\
=E[D_t^\bot(\sum_{k=0}^n\la L^{(k)},L^{(n-k)}\ra_T-\sum_{k=0}^n\la L^{(k)\bot},L^{(n-k)\bot}\ra_T)|\cF_t]\\
=2\sum_{k=0}^n[-(\a_k\a_{n-k}-\b_k\b_{n-k})W_t^\bot+(\a_k\b_{n-k}+\b_k\a_{n-k})W_t]\int_0^T(T-s)^{2n+2}ds\\
=\frac{2(T-t)^{2n+3}}{2n+3}\sum_{k=0}^n[-(\a_k\a_{n-k}-\b_k\b_{n-k})W_t^\bot+(\a_k\b_{n-k}+\b_k\a_{n-k})W_t].
\eeaa
Equalising coefficients at $W,W^\bot$ we obtain the desired formula.
One can be checked that $\lim_{n\to\infty}\sqrt[n]{|a_n|}=0,\lim_{n\to\infty}\sqrt[n]{|b_n|}=0.$
 Introducing $\a(s)=\sum_{n=0}^\infty\a_n  s^{2n+1},\;\b(s)=\sum_{n=0}^\infty\b_n  s^{2n+1}$ one obtains
 \beaa
 L_t=L_0+\int_0^t(\a(T-s)W_s+\b(T-s)W_s^\bot)dW_s,\\L_t^\bot=L_0^\bot+\int_0^t(\b(T-s)W_s-\a(T-s)W_s^\bot)dW_s^\bot.
 \eeaa
On the other hand we can derive ODE for the pair $(\a,\b)$
\bea\lbl{ode}
\a'(s)=2\a^2(s)-2\b^2(s),\;\a(0)=0,\\
\notag \b'(s)=1+4\a(s)\b(s),\;\b(0)=0.
\eea
Indeed
\beaa
\a'(s)=\a_0+\sum_{n=0}^\infty(2n+3)\a_{n+1}s^{2n+2}\\
=2\sum_{n=0}^\infty\sum_{k=0}^n(\a_k\a_{n-k}-\b_k\b_{n-k})s^{2n+2}=2a^2(s)-2\b^2(s),\\
\b'(s)=\b_0+\sum_{n=0}^\infty(2n+3)\b_{n+1}  s^{2n+2}\\
=1+4\sum_{n=0}^\infty\sum_{k=0}^n\a_k\b_{n-k}  s^{2n+2}=1+4\a(s)\b(s).
\eeaa
The equation  (\ref{ode}) is easy to solve, if we pass to the equation for complex-variable function $\z(s)=\a(s)+i\b(s)$
\beaa
\z'(s)=i+2\z^2(s), \;\z(0)=0.
\eeaa
 It is obvious that $\z(s)=\frac{1}{1-i}\tan((1+i)s)$ is a solution. We have
 \beaa\z(s)=\frac12(1+i)\frac{\sin((1+i)s)\cos((1-i)s)}{|\cos((1+i)s)|^2}
 \\=\frac14(1+i)\frac{\sin(2s)+i\sinh(2s)}{|\cos((1+i)s)|^2}\\=\frac14\frac{\sin(2s)-\sinh(2s)+i(\sin(2s)+\sinh(2s))}{\cos^2(s)\cosh^2(s)+\sin^2(s)\sinh^2(s)}.\eeaa
Finally we can  write  explicit solution
\beaa
\a(s)=\frac14\frac{\sin(2s)-\sinh(2s)}{\cos^2(s)\cosh^2(s)+\sin^2(s)\sinh^2(s)},\\
\b(s)=\frac14\frac{\sin(2s)+\sinh(2s)}{\cos^2(s)\cosh^2(s)+\sin^2(s)\sinh^2(s)}\eeaa
of (\ref{ode})  and conclude that it exists on whole $[0,\infty)$, since the denominator does not vanish.

\appendix

\section{Appendix}

 The formula $Ee^{-T^2\int_0^{1}W_t^2dt}=\frac1{\sqrt{\cosh(\sqrt2 T)}}$ is derived in \cite{ku}.
Similarly we can prove
\be{prop} $$Ee^{\int_0^{T}W_t^2dt}=\be{cases}\frac1{\sqrt{\cos(\sqrt2 T)}},\; {\rm if}\; T<\frac\pi{2\sqrt2}\\
\infty, \;{\rm if}\;T\ge\frac\pi{2\sqrt2}\ee{cases}.$$
\ee{prop}
{\it Proof.} Let $e_n(t)$ be orhonormal basis in $L^2[0,1]$. Then
$Ee^{\int_0^{T}W_t^2dt}=Ee^{T^2\int_0^{1}W_{t}^2dt}=Ee^{T^2\sum_{n=1}^\infty(\int_0^{1}e_n(t)W_{t}dt)^2}=E\prod_{n=1}^\infty e^{T^2(\int_0^{1}e_n(t)W_{t}dt)^2}$.
Since $$E(\int_0^{1}e_n(t)W_{t}dt)(\int_0^{1}e_m(t)W_{t}dt)=\int_0^Te_n(t)\int_0^T(t\wedge s)e_m(s)dsdt$$ it is convenient to use the orthonormal basis of eigenvectors of the operator $\int_0^T(t\wedge s)f(s)ds$ in $L^2[0,1]$.
From $\l f(t)=\int_0^T(t\wedge s)f(s)ds$ follows that $\l f''(t)=-f(t),\;f(0)=0,\;f'(1)=0$.
The function $\sin\mu\pi t$ satisfies  these conditions iff $\mu^2=1/\l$, $\cos\mu\pi=0$ and $\mu=-1/2+n$. Thus
$$\l_n=\frac1{(n-1/2)^2\pi^2},\;e_n(t)=\sqrt2\sin((n-1/2)\pi t),\;n\ge1$$ and $E(\int_0^{1}e_n(t)W_{t}dt)(\int_0^{1}e_m(t)W_{t}dt)=\l_n\int_0^1e_n(t)e_m(t)dt=0,\;n\neq m$. Since random variables
$(\int_0^{1}e_n(t)W_{t}dt)$ are orthogonal and normal they are  also independent.    Hence taking into account infinite product decomposition of $\cos(\sqrt2t)$ one gets
\beaa
Ee^{\int_0^{T}W_t^2dt}=\prod_{n=1}^\infty Ee^{T^2(\int_0^{1}e_n(t)W_{t}dt)^2}\\=\prod_{n=1}^\infty Ee^{T^2\l_nW_1^2}=\prod_{n=1}^\infty\frac1{\sqrt{1-\frac{2T^2}{(n-1/2)^2\pi^2}}}\\
=\sqrt{\prod_{n=1}^\infty\frac1{{1-\frac{8T^2}{(2n-1)^2\pi^2}}}}=\frac1{\sqrt{\cos(\sqrt2T)}},\eeaa
if $\sqrt2T<\pi/2$.

\noindent It easy to see that $$E\exp{\left(\int_0^{\frac{\pi}{2\sqrt2}}W_t^2dt\right)}=\lim_{T\uparrow\frac{\pi}{2\sqrt2}}E\exp{\left(\int_0^{T}W_t^2dt\right)}=\lim_{T\uparrow\frac{\pi}{2\sqrt2}}\frac1{\sqrt{\cos(\sqrt2T)}}=\infty.$$

If $T>\frac{\pi}{2\sqrt2}$ then $Ee^{\int_0^{T}W_t^2dt}>Ee^{\int_0^{\frac{\pi}{2\sqrt2}}W_t^2dt}=\infty.$

\be{lem} Let $(a_n)_{n\ge0}$ be a solution of the system
\bea\lbl{rec}
a_0=1,\;a_{n+1}=\sum_{k=0}^na_ka_{n-k}.\eea
 Then $a_n=\frac1{4n+2}\binom{2n+2}{n+1}$.
\ee{lem}
{\it Proof.}
For the series $u(\l)=\sum_{n=0}^\infty a_n\l^n$ from (\ref{rec}) we get equation $u(\l)=1+\l u^2(\l)$, with the roots $u(\l)=\frac1{2\l}(1\pm\sqrt{1-4\l})$.
The equality    $u(\l)=\frac1{2\l}(1+\sqrt{1-4\l})$ is impossible, since decomposition of the right  hand side is starting from the term $\frac1\l$.
Therefore, equality $a_n=\frac1{4n+2}\binom{2n+2}{n+1}$ follows   from the Taylor expansion of $1-\sqrt{1-4\l}$, since
\beaa
u(\l)=\frac1{2\l}(1-\sqrt{1-4\l})\\
=-\frac12\sum_{n\ge1}\frac{\frac12(\frac12-1)\cdots(\frac12-n+1)}{n!}(-4)^n\l^{n-1}\\
=\frac12\sum_{n\ge1}\frac{(2-1)\cdots(2n-2-1)}{2^nn!}
4^n\l^{n-1}\\
=\frac12\sum_{n\ge1}\frac{(2n-3)!!}{n!}2^{n}\l^{n-1}
=\frac12\sum_{n\ge1}\frac1{2n-1}\binom{2n}{n}\l^{n-1}.
\eeaa

\be{lem} There exist sequences $(m_i,i\ge1)\in\cM,\;(m_i^\bot,i\ge1)\in\cM^\bot,$ such that $e^\eta=c_1\frac{\cE_T(m_1)}{\cE_T(m_1^\bot)}\cE_T^2(m_1^\bot)$ and
\beq\lbl{fr}
e^\eta=c_n\frac{\cE_T(\sum_i^nm_i)}{\cE_T(\sum_i^nm_i^\bot)}\cE_T^2(m_n^{'\bot}),\;n\ge2,
\eeq
where $m_n^{'\bot}=m_n^\bot-\la m_n^\bot,\sum_i^{n-1}m_i^\bot\ra$.
\ee{lem}
{\it Proof.} The theorem will be proved by induction. Assume (\ref{fr}) is valid for n. There exist such martingales $m_{n+1},m_{n+1}^\bot$ that $c'\cE_T(m_{n+1}'+m_{n+1}'^\bot)=\cE_T^2(m_n^{'\bot})$
and
$$m_{n+1}'=m_{n+1}-\la m_{n+1},\sum_i^n,m_{i}\ra,\;\; m_{n+1}'^\bot=m_{n+1}^\bot-\la m_{n+1}^\bot,\sum_i^nm_{i}^\bot\ra$$
 are martingales w.r.t. $\cE(\sum_i^nm_i+m_i^\bot)\cdot P.$
Thus
\beaa
e^\eta=c_nc'\frac{\cE_T(\sum_i^nm_i)}{\cE_T(\sum_i^nm_i^\bot)}\cE(m_{n+1}-\la m_{n+1},\sum_i^nm_{i}\ra+m_{n+1}^\bot-\la m_{n+1}^\bot,\sum_i^nm_{i}^\bot\ra)\\
=c_{n+1}\frac{\cE_T(\sum_i^nm_i)\cE_T(m_{n+1}-\la m_{n+1},\sum_i^nm_{i}\ra)}{\cE_T(\sum_i^nm_i^\bot)\cE_T(m_{n+1}^\bot-\la m_{n+1}^\bot,\sum_i^nm_{i}^\bot\ra)}\cE_T^2(m_{n+1}^\bot-\la m_{n+1}^\bot,\sum_i^nm_{i}^\bot\ra)\\
=c_{n+1}\frac{\cE_T(\sum_i^{n+1}m_i)}{\cE_T(\sum_i^{n+1}m_i^\bot)}\cE_T^2(m_{n+1}^{'\bot}).
\eeaa
{\bf Remark.} If we will prove the convergence of series $\sum_im_i, \sum_im_i^\bot$, then $m_n^\bot\to0,m_n^{'\bot}\to0,\;\cE(m_n^{'\bot})\to1$ and $e^\eta=c\frac{\cE_T(\sum_i^\infty m_i)}{\cE_T(\sum_i^\infty m_i^\bot)}$.

\end{document}